\input amstex
\input pictex
\documentstyle{amsppt}
\magnification\magstep1

\define\vo{\text{\rm vol}}

\topmatter
\title
Equivalence of Geometric and Combinatorial\\
Dehn Functions
\endtitle
\author
Jos\'e Burillo and Jennifer Taback
\endauthor
\abstract
We prove that if a finitely presented group 
acts properly discontinuously, cocompactly and by isometries 
on a simply connected Riemannian manifold, then the Dehn function of the
group and the corresponding filling function of the manifold are
equivalent, in a sense described below. 
\endabstract
\address
Departament de Matem\`atiques, Universitat Aut\`onoma de Barcelona,
08193 Bellaterra, Spain
\endaddress
\email
burillo\@mat.uab.es
\endemail
\address
Dept. of Mathematics and Statistics, University at Albany, Albany, NY
12222
\endaddress
\email
jtaback\@math.albany.edu
\endemail
\endtopmatter

\document

\heading
1. Dehn functions and their equivalence
\endheading

Let $X$ be a simply connected 2-complex , and let $w$ be an
edge circuit in
$X^{(1)}$. If $D$ is a van Kampen diagram for $w$ (see \cite{5}), then
the area of $D$ is defined as the number of 2-cells on $D$, and the
area of $w$, denoted $a(w)$, is defined as the minimum of the areas of all van
Kampen diagrams for $w$. The Dehn function of $X$ is then defined to
be
$$
\delta_X(n)=\max a(w),
$$
where the maximum is taken over all loops $w$ of length $l(w)\le n$.

Given two functions $f$ and $g$ from $\Bbb{N}$ to $\Bbb{N}$ (or, more
generally, from $\Bbb{R}^+$ to $\Bbb{R}^+$), we say that $f\prec g$ if
there exist positive constants $A$, $B$, $C$, $D$, $E$ so that
$$
f(n)\le Ag(Bn+C)+Dn+E.
$$
Two such functions are called equivalent (denoted $f\equiv g$) if
$f\prec g$ and $g\prec f$. 
The Dehn function
is invariant under quasi-isometries: when one
considers the 1-skeleton of a complex as a metric space with the
path metric, where every edge has length one,
two complexes with quasi-isometric 1-skeleta have
equivalent Dehn functions (see \cite{1}).

Let $G$ be a finitely presented group, and let $\Cal{P}$ be a finite
presentation for $G$. Let $K=K(\Cal{P})$ be the 2-complex associated
to $\Cal{P}$, i.e. the 2-complex with a single vertex, an oriented
edge for every generator of $\Cal{P}$, and a 2-cell for every relator,
attached to the edges according to the spelling of the relator. Then
the Dehn function of $\Cal{P}$ is, by definition, the Dehn function
$\delta_{\tilde K}$ of the universal covering of $K$. 
Two finite presentations $\Cal{P}$ and $\Cal{Q}$ for the same group
$G$ 
yield 2-complexes $\widetilde{K}(\Cal{P})$ and
$\widetilde{K}(\Cal{Q})$ with quasi-isometric 1-skeleta,
and hence equivalent Dehn functions.  
Thus the Dehn function of the group $G$ is defined to be the
equivalence class of the Dehn function of any of its
presentations. An extensive treatment of Dehn functions of
finitely presented groups is given in \cite{4}.

A closely related definition can be formulated in the context of
Riemannian manifolds, dating back to the isoperimetric problem for
$\Bbb{R}^n$ in the calculus of variations. Given a Lipschitz loop $\gamma$
in a simply connected Riemannian manifold $M$, we define the area of
$\gamma$ to be the infimum of the areas of all Lipschitz discs bounded
by $\gamma$. We then define the geometric Dehn function of $M$ by
$$
\delta_M(x)=\max_{l(\gamma)\le x}\text{area}(\gamma)
$$  
where $l(\gamma) = length(\gamma)$.

It is natural to consider the question of whether the Dehn functions
of a simply connected Riemannian manifold $M$ and of a finitely presented
group $G$ acting properly discontinuously and cocompactly on $M$
agree. The fact that they effectively agree has been implicitly
assumed in the literature, though no proof has been given.  
A closely
related statement is given in \cite{2,~{\rm Theorem 10.3.3}}, 
applying the Deformation Theorem of Geometric Measure Theory (\cite{3,~4.2.9}
and \cite{7}) to this setting,
and which provides the basis of the Pushing Lemma below. This paper is
devoted to providing a complete and detailed proof that the
combinatorial and geometric Dehn functions
are equivalent. It is known to the authors that
M. Bridson has lectured on an alternate, 
unpublished proof for the same result. The
authors would like to thank Professor S. M. Gersten for his
encouragement and his useful remarks, Kevin Whyte for helpful
conversations and the referee for his precise
comments.

We prove the following theorem.  

\proclaim{Theorem 1.1} Let $M$ be a simply connected Riemannian
manifold, and $G$ a finitely presented group acting properly
discontinuously, cocompactly and by isometries on $M$. Let $\tau$ be a
$G$-invariant triangulation of $M$. Then the following three Dehn
functions are equivalent:
\roster
\item the Dehn function $\delta_G$ of any finite presentation of $G$,
\item the Dehn function $\delta_{\tau^{(2)}}$ of the 2-skeleton of
$\tau$, and
\item the geometric Dehn function $\delta_M$ of $M$.
\endroster
\endproclaim

The fact that $\delta_G$ and $\delta_{\tau^{(2)}}$ are equivalent is
clear: since $G$ acts cocompactly on $\tau$, there is a
quasi-isometry between $\tau^{(1)}$ and the 1-skeleton of
$\widetilde{K}(\Cal{P})$ for
any presentation $\Cal{P}$ of $G$, and the equivalence follows from
the results in \cite{1}. We will concentrate on the proof of the
equivalence between $\delta_{\tau^{(2)}}$ and $\delta_M$.  
The arguments will be mainly geometric, relating the lengths and
areas of loops and discs in $M$ with those included in the
triangulation $\tau$. The first step in this direction is the
Pushing Lemma, which is a complete analog of the Deformation Theorem in
Geometric Measure Theory and already stated and proved, in a slightly
different way, in \cite{2,~{\rm Theorem 10.3.3}}, and whose proof we will 
follow closely.

\heading
2. Technical Lemmas
\endheading

The Pushing Lemma, stated below, will allow us to relate arbitrary Lipschitz
chains in $M$ to chains in the corresponding skeleta of $\tau$. The main
technical problem to be overcome is that projection of a Lipschitz
chain to $\tau$ from a badly
chosen point can increase the volume of the chain arbitrarily. We overcome
this by using techniques from measure theory that assure the existence of a
center of projection far enough from the chain, thus
providing control on the growth of the volume.

\vfil
\pagebreak

\proclaim{Lemma 2.1 (Pushing Lemma)} Let $M$, $G$ and $\tau$ be as
above. Then there exists a constant $C$, depending only on $M$ and 
$\tau$, with the following property:
Let $T$ be a Lipschitz $k$-chain in $M$, such that $\partial T$ is
included in $\tau^{(k-1)}$. Then there exists another Lipschitz
$k$-chain $R$, with $\partial R=\partial T$, which is included in
$\tau^{(k)}$, and a Lipschitz $(k+1)$-chain $S$, with $\partial
S=T-R$, satisfying
$$
\vo_k(R)\le C\vo_k(T)\qquad \text{and}\qquad
\vo_{k+1}(S)\le C\vo_k(T).
$$
In particular, if $T$ is a loop, so is $R$, and $S$ is a homotopy
from $T$ to $R$.
\endproclaim

The Pushing Lemma differs from the statement in \cite{2} because it
applies to chains as well as cycles, since
the boundary of the chain is not modified, as it is included in the
$(k-1)$-skeleton.  
A statement for
cycles is not sufficient, since this lemma will be applied to
chains as well as loops, and the fact that $\partial T=\partial
R$ is crucial in the proof of the main theorem.

We first prove a lemma which will later allow us to choose our center of
projection to lie away from the Lipschitz chain $T$.

\proclaim{Lemma 2.2}
Let $f: S^k \rightarrow \sigma_{k+1}$ be Lipschitz with constant $L$,
where $\sigma_{k+1}$ is the standard Euclidean $(k+1)$-simplex.  Then
$f(S^k)$ has Lebesgue $(k+1)$-measure zero.
\endproclaim

\demo{Proof}
Since $S^k$ is compact, choose a finite open cover of $S^k$ by
$k$-dimensional balls $B_i$
of radius $\frac{1}{n}$.  We can cover $S^k$ with  $C_1n^k$ such balls, for some constant $C_1$.  
The image of any ball $B_i$ under the Lipschitz map $f$ is contained in a
$(k+1)$-dimensional ball $B_i' \subset \sigma_{k+1}$ with
$(k+1)$-volume $\frac{C_2}{n^{k+1}}$ for some constant $C_2$.  
Then the total volume of the collection $\{ B_i' \}$ is at most
$\frac{C_1C_2}{n}$.  
So $f(S^k)$ is contained in an open set of $\sigma_{k+1}$ whose total
volume is $\frac{C_1C_2}{n}$ and thus $f(S^k)$ has Lebesgue measure
$0$.  
\enddemo

\demo{Proof of Lemma 2.1}
The proof will proceed by descending induction on the skeleta of
$\tau$. Assume that a Lipschitz $k$-chain $T$ is included in $\tau^{(i)}$ but not in
$\tau^{(i-1)}$, for $i>k$. We want to proceed simplex by simplex,
choosing an appropriate point not in $T$ in each simplex and projecting  the
chain $T$ radially
from this point to the boundary of the simplex. 
We will prove the following claim.

\medskip
\noindent
{\it Claim:} 
There exists a constant $C$ with the property that for every simplex,
there is a point $p$ not in $T$ so that radial projection of $T$ 
from $p$ to the
boundary of the simplex does not increase the volume of the chain by
more than a multiplicative factor $C$. 

\medskip

Observe that since
$T$ is compact, it only intersects finitely many simplices of $\tau$, 
and in each simplex is only modified by a radial projection from
a point not in $T$. These radial projections only increase the
Lipschitz constant of $T$, but the chain $R$ obtained after the 
projections will still be Lipschitz.

To simplify the computations, we will work through the proof in the unit 
Euclidean simplex of edge length one. Since $G$ acts cocompactly on $M$, 
we can construct a sufficiently fine finite triangulation of the quotient and
lift it to $M$. If the simplices are small enough we can map them to 
$\Bbb{R}^n$ via the exponential map.  
Since the 
exponential map is Lipschitz, the changes in the metric are bounded
by only a multiplicative constant. We then have a finite number of simplices
in $\Bbb{R}^n$, so the distortion is again bounded. Thus working with the
unit simplex only affects the value of the constant $C$.

Let $\sigma$ be the unit Euclidean $i$-simplex, 
$O$ the barycenter of $\sigma$, and
$r$ a positive number so that the ball of center $O$ and
radius $3r$ is included in the interior of $\sigma$. Let $B$ be the
ball of center $O$ and radius $r$, with $u$ an element of $B$, and $B_u$ the
ball of center $u$ and radius $2r$. Clearly $B\subset B_u$, for all
$u$. Let $\pi_u$ be the radial projection with center $u$ of
$B_u\setminus \{u\}$ onto $\partial B_u$. Let $Q=T\cap\sigma$.
We want to see that there exists a constant $v_0$ independent of $T$
and $\sigma$, and a point $u \in B \setminus Q$, dependent on $T$,
with
$$
\vo_k(\pi_uQ)\le v_0\,\vo_k(Q).
$$
From Lemma 2.2, we see that the set $B \setminus Q$ has the same
measure as $B$, allowing us to choose $u \in B \setminus Q$.  

For every positive real number $v$ define
$$
A_v=\{u\in B \setminus Q\,|\,\vo_k(\pi_uQ)>v\,\vo_k(Q)\}
$$
and let $\alpha(v)=m_i(A_v)$, where $m_i$ is the $i$-dimensional
Lebesgue measure. We want to prove that
$$
\lim_{v\to\infty}\alpha(v)=0.
$$
Then we will choose $v_0$ with $\alpha(v_0)<m_i(B)$, so the measure of
$A_{v_0}$ will be less than the measure of $B$.  
Thus there will
exist a point $u\in (B\setminus Q) \setminus A_{v_0}$,
which will be the center of projection. Since
$u\notin A_{v_0}$, this projection will increase the area at most by a
multiplicative factor $v_0$.

\beginpicture
\setcoordinatesystem point at -181 110
\setplotarea x from -181 to 100, y from -85 to 100
\circulararc 360 degrees from 15 0 center at 0 0
\circulararc 240 degrees from -7 -37 center at -7 -7
\setlinear 
\plot 0 100  87 -50  -87 -50  0 100 /
\setquadratic
\setdashes <3pt>
\plot -33 8  -15 -15  -7 -37 /
\setsolid
\plot -50 14  -40 13  -33 8 /
\plot -7 -37  -5 -43  0 -50 /
\circulararc 120 degrees from -33 8 center at -7 -7
\put {$\sigma$} [lB] at 18 75
\put {$O$} [lB] at 1 1
\put {$u$} [l] at -5 -7
\put {$B$} [lt] at 0 -17
\put {$B_u$} [lt] at 15 -30
\put {$Q$} [rt] at -16 -16
\put {$\pi_u Q$} [rt] at -30 -30
\put {.} at 0 0
\put {.} at -7 -7
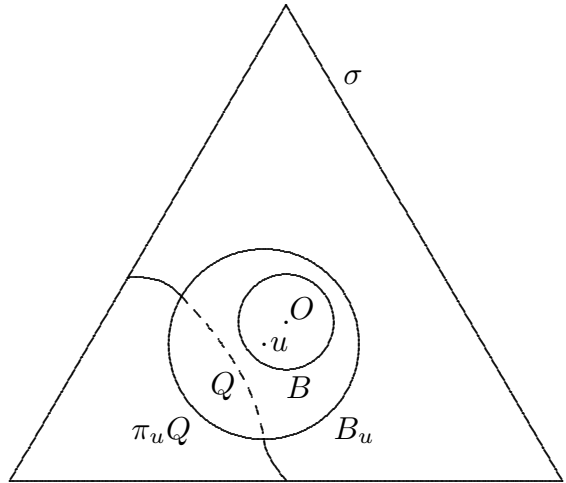
\put {Figure 1: \sl Projecting $Q$ to the boundary of $B_u$.} [B] at
0 -70
\endpicture

We have
$$
\split
\vo_k(\pi_uQ)&\le\vo_k(\pi_u(Q\cap B_u))+\vo_k(Q)\\
&\le\int_{Q\cap B_u}\left(\frac{2r}{||x-u||}\right)^k\,dx+\vo_k(Q),
\endsplit
$$
where the first term accounts for the volume obtained after
projecting, and the second term takes care of the possibility of $Q$
and $B_u$ being disjoint. Assume now that $\vo_k(Q)$ is nonzero (if
$\vo_k(Q)=0$ then $\vo_k(\pi_uQ)=0$). Then we have:
$$
\split
\alpha(v)\,v\,\vo_k(Q)&=v\,\vo_k(Q)\int_{A_v}du=\int_{A_v}v\,
\vo_k(Q)\,du\\
&\le\int_{A_v}\vo_k(\pi_uQ)\,du\le\int_B\vo_k(\pi_uQ)\,du\\
&\le\int_B\left(\int_{Q\cap B_u}\left(\frac{2r}{||x-u||}\right)^k\,dx
+\vo_k(Q)\right)\,du\\
&=(2r)^k\int_{Q\cap B_u}\int_B ||u-x||^{-k}\,du\,dx+
\vo_i(B)\vo_k(Q).
\endsplit
$$
Notice that the function $\left(\frac{2r}{||x-u||}\right)^k$ 
is bounded above and
below, since $u \notin Q \cap B_u$, and is integrated over compact
regions.  This allows us to change the order of integration.  
Now make a change of variables, letting $w=u-x$, and increase the domain
of integration to $B(0,3r)$.  
We continue with the upper bound for $\alpha(v)\,v\,\vo_k(Q)$:
$$
\split
\alpha(v)\,v\,\vo_k(Q)&\le (2r)^k\int_{Q\cap B_u}\int_B ||u-x||^{-k}\,du\,dx+
\vo_i(B)\vo_k(Q)\\
&\le(2r)^k\int_{Q\cap B_u}dx\int_{B(O,3r)}||w||^{-k}\,dw+
\vo_i(B)\vo_k(Q)\\
&\le K\vo_k(Q),
\endsplit
$$
where
$$
K=(2r)^k\int_{B(O,3r)}||w||^{-k}\,dw+\vo_i(B).
$$
Observe that $K$ is finite and independent of $T$ and $\sigma$. 
We conclude that $\alpha(v)v\le K$. Knowing $K$, we can find $v_0$
such that $K/v_0<m_i(B)$, where $v_0$ is a constant
independent of $T$ and $\sigma$. We have now found $A_{v_0}$ with
strictly less measure than $B$, and can pick a point in $(B\setminus
Q) \setminus
A_{v_0}$ from which to project so that the volume increases
at most by a multiplicative factor $v_0$.

The result of the above argument is the construction of another chain
$\pi_uQ$ which is far enough from $O$. We can now project radially
from $O$ to $\partial\sigma$, and the change of volume is bounded
since $\pi_uQ$ is at least at a distance $r$ from $O$. The combination
of this change of volume with $v_0$ gives the constant needed in this
precise skeleton. Combining the constants from all of these steps, we
obtain the desired constant $C$.
Observe that these projections leave $\tau^{(i-1)}$ unchanged, so
clearly $\partial T$ is preserved.

The $(k+1)$-chain $S$ is obtained by joining every $x\in Q$ to
$\pi_ux$ by a segment. The volume of the piece of $S$ contained in
$\sigma$ is then bounded, as before, by
$$
(2r)^{k+1}\int_{Q\cap B_u}\frac{dx}{||x-u||^k},
$$
where the extra factor $2r$ is obtained from the direction of the
projection, since each segment has length bounded by $2r$. An argument
similar to the previous one shows that projecting from most points
in $B$ gives the right bound for the volume.  $\square$
\enddemo

The third lemma states that for a Lipschitz map, almost every point
in the target space has a finite number of preimages. It is a direct
consequence of the area formula for Lipschitz maps, and it will be
used for both loops and discs in the proof of Theorem 1.1.

\proclaim{Lemma 2.3} Let $T$ be a Lipschitz $k$-chain in $M$, where
$k\le \dim M$. Then the set of points in $M$ with infinite preimages
under $T$ has Hausdorff $k$-measure zero.
\endproclaim

\demo{Proof}
Let $\sigma_k$ be the standard closed $k$-simplex, and let 
$$
E:\sigma_k\longrightarrow M
$$
be one of the simplices in $T$. Since $E$ is a Lipschitz map, 
by Rademacher's Theorem (\cite{3,~3.1.6}) it is
differentiable almost everywhere (with respect to the Lebesgue
$k$-measure), so the Jacobian $J_k E(x)$ is well defined for almost
all $x\in \sigma_k$. For $y\in M$,
let $N(E,y)$ be the number of elements of
$E^{-1}(y)$, possibly infinite,
and denote by $m_k$ and $h_k$ the Lebesgue and Hausdorff $k$-measures, respectively. Then the area formula for Lipschitz maps
(\cite{3,~3.2.3}) states that
$$
\int_{\sigma_k} |J_k E(x)|\,dm_k(x)=\int_M N(E,y)\,dh_k(y).
$$
Since $E$ is Lipschitz, we know that $|J_k E(x)|$ is bounded,
and since $\sigma_k$ has finite measure, the
integral on the left hand side is finite. So the set where $N(E,y)$ is
infinite cannot have positive Hausdorff $k$-measure, because then the right
hand side of the equation would be infinite. $\square$
\enddemo

\heading
3. Proof of the Main Theorem
\endheading

We begin by proving the one of the two inequalities necessary for the equivalence of $\delta_M$ and
$\delta_{\tau^{(2)}}$, namely
$$
\delta_M\prec\delta_{\tau^{(2)}}.\tag3.1
$$
Let $\gamma$ be a Lipschitz loop in $M$, with length at most
$n$. Using the Pushing Lemma, we can construct a new loop $\eta$, of
length at most $Cn$, which is included in the 1-skeleton,
and the homotopy between $\gamma$ and $\eta$ has
area at most $Cn$.

The loop $\eta$ is not necessarily combinatorial, but it is a rectifiable loop
in a non-positively curved space, namely the metric graph $\tau^{(1)}$. So
there is a unique (up to reparametrization) closed geodesic $\zeta$ in the free
homotopy class of $\eta$. The straight homotopy (in $\tau^{(1)}$)
from $\eta$ to $\zeta$ is a
map from an annulus to $\tau^{(1)}$. The length of $\zeta$ decreases
monotonically and its area can be made arbitrarily small.

The combinatorial loop $\zeta$ can be
filled combinatorially by at most $\delta_{\tau^{(2)}}(Cn)$
2-simplices in $\tau$. Thus
$$
\delta_M(n)\le A\delta_{\tau^{(2)}}(Cn)+2Cn,
$$
where $A$ is the area of the largest 2-cell in $\tau$, and 
it follows that $\delta_M\prec\delta_{\tau^{(2)}}$.

To prove the reverse inequality
$$\delta_{\tau^{(2)}}\prec\delta_M,$$
to (3.1), we start with a combinatorial loop 
$\gamma$ in the 1-skeleton of $\tau$, with length at most $n$. Let
$$
f:D^2\longrightarrow M
$$
be a Lipschitz disc in $M$ with boundary $\gamma$, and with area $a$. We want to
construct a van Kampen diagram for $\gamma$ and bound its area in terms
of $a$. The first step is to use the Pushing Lemma to find
a new disc (also denoted $f$) which is included in $\tau^{(2)}$, and
whose area is at most $Ca$.

Let $\sigma$ be an open 2-simplex of $\tau$ contained in $f(D^2)$. By Lemma 2.3,
we can choose a point $p\in\sigma$, such that $f^{-1}(p)$ is
finite. Let $X$ be a component of $f^{-1}(\sigma)$. 
If $X\cap f^{-1}(p)=\varnothing$, then $f\big|_{\dsize X}$ can be
modified by composing with a radial projection from $p$. After this
change, a component $X$ of $f^{-1}(\sigma)$ satisfies $X\cap 
f^{-1}(p)\ne\varnothing$, and there are only finitely many of these
components. Moreover, if $f\big|_{\dsize X}$ is not surjective, we
can again modify $f\big|_{\dsize X}$ by a radial projection from a point
not in $f(X)$, to push its image to $\partial\sigma$. After these
changes to $f$, there is a component $X$ of $f^{-1}(\sigma)$ so that
$f\big|_{\dsize X}$ is surjective, and $X\cap f^{-1}(p)\ne\varnothing$.
If $X$ is one such component, the original $f$ has not been modified
in $X$ by any radial projection, and the map
$$
f\big|_{\dsize X}:X\longrightarrow\sigma
$$
is still Lipschitz, since it is the restriction of the original map $f$.

We will obtain a lower bound on the area of $f\big|_{\dsize X}$ using the degree of $f\big|_{\dsize X}$.
Since $f\big|_{\dsize X}$ is differentiable almost everywhere,
we can define the degree of $f\big|_{\dsize X}$ at a
point $y\in f(X)$ by
$$
\text{deg}\,f\big|_{\dsize X} (y) = \sum_{x\in f^{-1}(y)}
\text{sign}\,J_2f(x).
$$
Moreover, since $X$ is an open  
connected component of $f^{-1}(\sigma)$, we
have that $f(X)\subset\sigma$ and $f(\partial X)\subset\partial\sigma$, so
$f(X)$ and $f(\partial X)$ are disjoint. Then, by \cite{3,~4.1.26}, 
the degree of $f\big|_{\dsize X}$ 
is almost constant in $f(X)$, and we can define the degree of
$f\big|_{\dsize X}$ as the value $d_X$ it achieves at almost every
$y\in f(X)$. The lower bound on the area of $f\big|_{\dsize X}$ is given by the area formula for
Lipschitz maps: if $u$ is an integrable function with respect to $m_2$, we
have (see \cite{3,~3.2.3}):
$$
\int_X u(x)|J_2f(x)|\,dm_2=\int_{\sigma} \sum_{x\in f^{-1}(y)\cap X} 
u(x)\,dh_2,
$$
and taking $u(x)=\text{sign}\,Jf(x)$ we obtain:
$$
\gather
\text{area}\, f\big|_{\dsize X} =\int_X |J_2f(x)|\,dm_2\ge\left|\int_X
J_2f(x)\,dm_2\right|=\\
\left|\int_X \text{sign}\,J_2f(x)\,|J_2f(x)|\,dm_2\right|=
\left|\int_{\sigma}
\text{deg}\,f\big|_{\dsize X}\,dh_2\right|=\frac{\sqrt 3}4 |d_X|.
\endgather
$$

Our goal is to find a simplicial map
$$
g:D^2\longrightarrow \tau^{(2)}
$$
(with some simplicial structure in $D^2$) such that only $|d_X|$ simplices are mapped by the identity to $\sigma$
under $g\big|_{\dsize
X}$, and the
rest of $X$ is mapped to $\partial\sigma$. Then we will have that the
combinatorial area of $g$ is bounded as follows,
$$
\sum_X |d_X|\le \sum_X \frac4{\sqrt 3}\text{area}\left(f\big|_{\dsize
X}\right)\le\frac 4{\sqrt 3}Ca
$$
giving us the required bound. Note that the map $g$
is not combinatorial, but only simplicial, and at the end of the proof a short
argument will be required to ensure the existence of a combinatorial map
whose area admits the same upper bound.

The first step in finding the map $g$ is to smooth the map 
$f\big|_{\dsize X}$, in order to apply differentiable techniques to it. Let $O$ be the barycenter
of $\sigma$, and choose $0<\epsilon<r$ such that:
$$
\varnothing\ne B(O,r-\epsilon)\subset B(O,r)\subset B(O,2r)
\subset B(O,2r+\epsilon)\subset\sigma,
$$
and let $U_1=f^{-1}(B(O,r))$ and $U_2=f^{-1}(B(O,2r))$. We have that
$\overline{U_1}\subset U_2\subset\overline{U_2}\subset X$. Choose
$\delta>0$ so that $B(x,\delta)\subset X$ for all $x\in U_2$, and
so that if $|x-y|<\delta$ then $|f(x)-f(y)|<\epsilon$, for all
$x,y\in X$. Let $\varphi$ be a $C^\infty$ bump function in $\Bbb{R}^2$
with support in $B(0,\delta)$, and with integral 1. Then, 
for $x\in U_2$, we can construct the convolution
$$
f*\varphi(x)=\int_{B(x,\delta)}f(x-z)\varphi(z)\,dz,
$$
which is $C^\infty$ in $U_2$, and satisfies
$|f(x)-f*\varphi(x)|<\epsilon$ for all $x\in U_2$.
Also, if $f\big|_{\dsize X}$ 
was Lipschitz with constant $L$, then $f*\varphi$ is also
Lipschitz with the same constant: if $x,y\in U_2$,
$$
|f*\varphi(x)-f*\varphi(y)|\le|f(x-z)-f(y-z)|\int_{B(0,\delta)}
\varphi(z)\,dz\le L|x-y|.
$$
Now choose a Lipschitz function $\alpha$ on $X$ with values in 
$[0,1]$ which is equal to 1 in $U_1$ and equal to 0 outside $U_2$, and define
$$
\tilde f =\alpha (f*\varphi) + (1-\alpha)f\big|_{\dsize X}.
$$
Note that $\tilde f$ is defined only on $X$.
Then $\tilde f$ satisfies the following properties:
\roster
\item $|f(x)-\tilde f(x)|<\epsilon$ for all $x\in X$,
\item $\tilde f$ is smooth in $U_1$,
\item $\tilde f=f$ in $X\setminus U_2$,
\item $\tilde f$ is Lipschitz, and
\item $\text{deg}\,\tilde f=\text{deg}\,f\big|_{\dsize X}$.
\endroster
The first three properties are clear from the construction of $\tilde f$, and
property (4) holds because 
$f\big|_{\dsize X}$ and $f*\varphi$ and $\alpha$ are all
Lipschitz. To see that the degree is unchanged, since the degree is
almost constant, and $f\big|_{\dsize X}$
and $\tilde f$ agree outside $U_2$, we only
need to find a point in $\sigma\setminus B(O,2r+\epsilon)$ for which
the degree is $d_X$ for both $f\big|_{\dsize X}$ and $\tilde f$.

We can now use Sard's Theorem (\cite{6})
to claim the existence of a regular
value for $\tilde f$ in $B(O,r-\epsilon)$ whose preimages are all in
$U_1$. Let $q$ be this regular value and let $p_1,\ldots,p_m$ be its
preimages. Let $V$ be an open disc with center $q$ such that
$\tilde f^{-1}(V)=V_1\cup\ldots\cup V_m$, where the $V_i$ are discs
around $p_i$, pairwise disjoint,
and such that $\tilde f\big|_{\dsize V_i}$ is a diffeomorphism. In
general, we will have that $m>|d_X|$, and must cancel
discs with opposite orientations. Assume $V_{m-1}$ and $V_m$ are
mapped to $V$ with opposite orientations. Choose $a\in\partial
V_{m-1}$ and $a'\in\partial V_m$ with $\tilde f(a)=\tilde f(a')$, and
join $a$ and $a'$ with a simple path $\lambda$ such that
$\tilde f(\lambda)$ is nullhomotopic in $\sigma\setminus V$.  
This can be done because the map
$$
\tilde f:X\setminus \bigcup_{i=1}^{m}V_i\longrightarrow
\sigma\setminus V
$$
induces a surjective homomorphism of fundamental groups.
After contracting $\tilde f(\lambda)$, we can assume $\tilde f(\lambda)$
is the constant path $\tilde f(a)$. Remove the discs
$V_{m-1}$ and $V_m$ and perform surgery along $\lambda$. The new
boundary thus created is mapped to $\partial V$ under $\tilde f$ by a
map from $S^1$ to itself of degree zero. Extend this map to a map from
$D^2$ to $S^1$ and attach it to $\tilde f$ along this boundary. For the
new map (which we will continue calling $\tilde f$),
the preimage of $q$ consists only of the points
$p_1,\ldots,p_{m-2}$. Repeating this process we will obtain a map
where now only the discs $V_1,\ldots,V_{|d_X|}$ are mapped to
$V$, all with the same orientation.

Choose (temporarily) a sufficiently fine subdivision of $\tau$ so
that there is a 2-simplex $W$ in $V$, and let
$\rho_i={\tilde f}^{-1}(W)$.
Modify the map in $X$ by composing with the expansion of $W$ into all
of $\sigma$.

\beginpicture
\setcoordinatesystem point at 0 78
\setplotarea x from 0 to 361, y from -103 to 70
\setdots <3pt>
\setlinear
\plot 40 68  160 68 /
\plot 30 51  170 51 /
\plot 20 34  180 34 /
\plot 10 17  190 17 /
\plot 0 0  180 0 /
\plot 10 -17  170 -17 /
\plot 20 -34  120 -34 /
\plot 30 -51  110 -51 /
\plot 40 -68  100 -68 /

\plot 0 0  40 68 /
\plot 10 -17  60 68 /
\plot 20 -34  80 68 /
\plot 30 -51  100 68 /
\plot 40 -68  120 68 /
\plot 60 -68  140 68 /
\plot 80 -68  160 68 /
\plot 100 -68  170 51 /
\plot 150 -17  180 34 /
\plot 170 -17  190 17 /

\plot 0 0  40 -68 /
\plot 10 17  60 -68 /
\plot 20 34  80 -68 /
\plot 30 51  100 -68 /
\plot 40 68  110 -51 /
\plot 60 68  120 -34 /
\plot 80 68  130 -17 /
\plot 100 68  150 -17 /
\plot 120 68  170 -17 /
\plot 140 68  180 0 /
\plot 160 68  190 17 /

\setsolid
\setquadratic
\plot 100 60  40 40  30 0  70 -60  100 -30  110 -10  130 0  160 10
      150 30  120 54  100 60 /

\setlinear
\plot 249 10  243 -5  255 -5  249 10 /
\plot 249 30  223 -15  275 -15  249 30 /
\plot 334 30  308 -15  360 -15  334 30 /

\arrow <7pt> [.25,.5] from 190 0 to 223 0
\arrow <7pt> [.25,.5] from 275 0 to 308 0

\plot 60 34  50 17  70 17  60 34 /
\plot 90 17  100 0  110 17  90 17 /
\plot 70 -17  60 -34  80 -34  70 -17 /

\put {$X$} [lB] at 139 44
\put {$\rho_1$} at 60 22.5
\put {$\rho_2$} at 100.5 11.5
\put {$\rho_3$} at 70 -28.5
\put {$W$} [lt] at 255 0
\put {$\sigma$} [lB] at 265 10
\put {$\sigma$} [lB] at 350 10
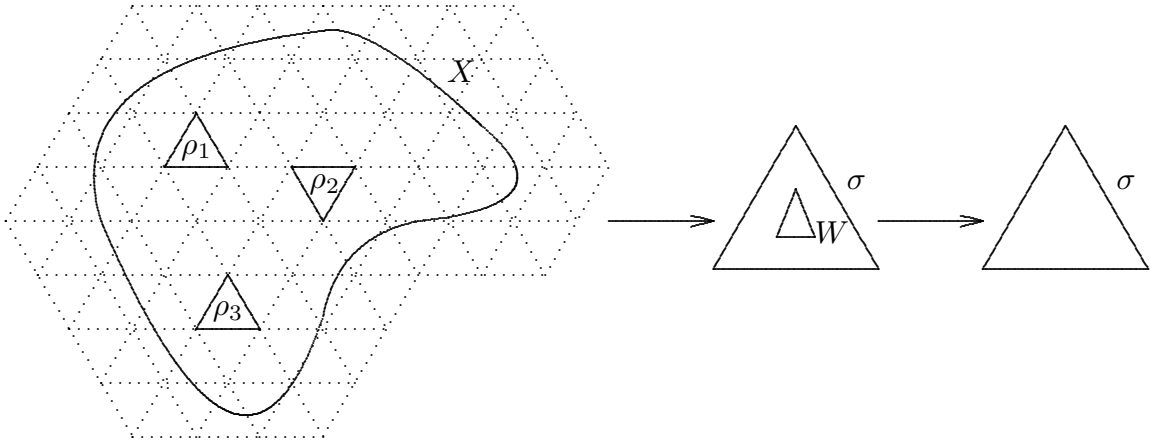
\put {Figure 2: \sl Making the map $f$ simplicial} [B] at 180.5 -88

\endpicture

After this process is done 
for all $\sigma$, we obtain a map from $D^2$ to $\tau^{(2)}$, where
all the $\rho_i$ are sent homeomorphically to
2-simplices of $\tau$, and the rest of $D^2$ is sent to the 1-skeleton
of $\tau$.
To finish the construction
of $g$, find a simplicial structure on $D^2$ compatible with the
simplicial structure on the original loop $\gamma$ and which includes
all the $\rho_i$ obtained for all $\sigma$ as 2-simplices.  
Now approximate the map $\tilde f$ simplicially within $\tau^{(1)}$ relative to all the $\rho_i$ and to $\gamma$.
The result is simplicial, and the number of simplices sent by $g$
homeomorphically to 2-simplices in $\tau$ is
$$
\sum_X|d_X|\le\frac4{\sqrt 3}Ca.
$$
This map
is not a van Kampen diagram yet, since it is only simplicial. To finalize the
proof of the inequality 
$$
\delta_{\tau^{(2)}}\prec\delta_M,
$$
we will find a van Kampen diagram which satisfies the same upper bound as
the map $g$. Consider simplicial
maps from a contractible planar 2-complex $Y$ into $\tau^{(2)}$,
with boundary $\gamma$, whose area satisfies the same bound as $g$.
(The map $g$ shows the existence of such maps.) 
Among all these maps, choose one with the minimum
number of 2-cells in $Y$. 
This map is necessarily combinatorial, since if some
2-cell of $Y$ is collapsed to the 1-skeleton of $\tau^{(2)}$, we could collapse it in $Y$
and find a map with fewer 2-cells. This map is the required van Kampen diagram
for the loop $\gamma$, and the second inequality is proved.

\enddocument